\documentclass[a4paper,Haag duality]{mathscan}
\newtheorem{thm}{Theorem}[section] 
\newtheorem{pro}[thm]{Proposition}

\theoremstyle{definition}           
\newtheorem{rem}[thm]{Remark}       

\newcommand{\NI}{\noindent}

\newcommand{\bea}{\begin{eqnarray}}
\newcommand{\eea}{\end{eqnarray}}
\newcommand{\dsp}{\displaystyle}
\def \b #1 {\bf #1}
\newcommand{\IR}{I\!\!R}
\newcommand{\IN}{I\!\!N}

\newcommand{\IC}{I\!\!C}

\newcommand{\IT}{I\!\!T}
\newcommand{\IZ}{Z\!\!\!Z}
\newcommand{\cal}{\mathcal}

\newcommand{\clm}{{\cal M}}

\newcommand{\clf}{{\cal F}}
\newcommand{\clg}{{\cal G}}
\newcommand{\clh}{{\cal H}}
\newcommand{\clp}{{\cal P}}

\newcommand{\clb}{{\cal B}}

\newcommand{\clj}{{\cal J}}
\newcommand{\cln}{{\cal N}}

\newcommand{\al}{\alpha}

\newcommand{\raro}{\rightarrow}

\newcommand{\vsp}{\vskip 1em}

\newcommand{\ul}{\underline}

\def \qed {\hfill \vrule height6pt width 6pt depth 0pt}
\newcommand{\be}{\begin{equation}}
\newcommand{\ee}{\end{equation}}
\newcommand{\ben}{\begin{eqnarray*}}
\newcommand{\een}{\end{eqnarray*}}
\begin{document}

\title{ Pure inductive limit state and Kolmogorov's property-II }

\author{ Anilesh Mohari }
\thanks{...}

\address{ The Institute of Mathematical Sciences, CIT Campus, Taramani, Chennai-600113 }

\email{anilesh@imsc.res.in}



\thanks{...}

\begin{abstract}
A translation invariant state $\omega$ on $C^*$-algebra $\clb=\otimes_{k \in \IZ}\!M^{(k)}$, where $\!M^{(k)}=\!M_d(\IC)$ is the $d-$dimensional matrices over field of complex numbers, 
give rises a stationary quantum Markov chain and associates canonically a unital completely positive normal map $\tau$ on a von-Neumann algebra $\clm$ with a faithful normal invariant 
state $\phi$. We give an asymptotic criteria on the Markov map $(\clm,\tau,\phi)$ for purity of $\omega$. Such a pure $\omega$ gives only type-I or type-III factor $\omega_R$ once restricted to 
one side of the chain $\clb_R=\otimes_{\IZ_+}\!M^{(k)}$. In case $\omega_R$ is type-I, $\omega$ admits Kolmogorov's property.  
\end{abstract}

\maketitle 

\section{ Introduction }

\vsp 
Let $\clb=\otimes_{\IZ} \!M_d(\IC)$ be the uniformly hyper-finite $C^*$-algebra over the lattice $\IZ$, where $\!M_d(\IC)$ be the $d \times d$-matrices over complex field $\IC$. 
A state $\omega$ on $\clb$ is called translation invariant if $\omega(x)=\omega(\theta(x))$ where $\theta$ is the translation induced by $z \raro  z + 1$ for all $z \in \IZ$. A state 
$\omega$ on $\clb$ is called a factor state if $\pi_{\omega}(\clb)''$ is a factor i.e. it's center is trivial, where $(\clh,\pi,\Omega)$ is the GNS space associated with $\omega$ on $\clb$ [BR1]. It is well known since late 60's [Pow] that a translation invariant state $\omega$ 
on $\clb$ is a factor state if and only if 
\be 
\mbox{sup}_{x \in \clb_{\Lambda_n^c}, ||x|| \le 1}|\omega(xy)) - \omega(x)\omega(y)| \raro 0
\ee 
for all $y \in \clb$ as $n \raro \infty$, where $\Lambda_{n}$ is the local algebra with support in the finite set $\{m: -n \le m \le n \}$. Further this criteria is equivalent to factor property of the state $\omega_R$, the restriction of the state $\omega$ to $\clb_R=\otimes_{\mathbb{Z}}\!M_d(\IC)$. Such an elegant criteria on $(\clb_R,\theta_R,\omega_R)$ were missing for purity of $\omega$. It is also known [BJKW,Ma] that for a type-I factor state $\omega_R$, $\omega$ is pure. There are pure $\omega$ where $\omega_R$ is a type-III factor state [AMa].     

\vsp 
Let $(\clh,\pi,\Omega)$ be the GNS space of $(\clb_R,\omega_R)$. We set support projection $P_0=[\pi(\clb_R)'\Omega]$ for the state $\omega_R$ in von-Neumann algebra $\pi(\clb_R)''$. 
Then by invariance property of the state $\omega_R(P_0\theta_R(I-P_0)P_0)=0$ which says that $\omega_R(P_0) \ge P_0$. We set unital completely positive map $\tau:\clm \raro \clm$ defined by 
\be
\tau(x)=P_0\theta_R(P_0xP_0)P_0
\ee   
for $x \in \clm_0=P_0\pi_{\omega_R}(\clb_R)''P_0$ with faithful normal invariant state 
$\phi_0(x)=<\Omega,x\Omega>$. Here we aim to continue line of investigations [Mo1,Mo2,Mo3] to find a useful necessary and sufficient condition in terms of asymptotic relation on $(\clm_0,\tau,\phi_0)$ for $\omega$ to be pure. In section 2 we deal with stationary Markov process associated with a Markov map $(\clm,\tau,\phi)$ and prove basic technical results in Theorem 2.1 and Theorem 2.4 which gives asymptotic criteria on $(\clm_0,\tau,\phi_0)$ for purity of the stationary state. We also deal with time reversed  
stationary Markov process which helped to find exact relation between purity of the stationary state and Power's shift criteria [Po2]. Theorem 2.5 gives a surprising result which says that pure stationary states once restricted to forward algebra, gives either a type-I or a type-III factor. In section 3 besides finding a criteria for purity of $\omega$ on $\clb$ in a more general $C^*$-algebraic set of inductive limit state, we prove the following theorem. 

\begin{thm} 
Let $\omega$ be a translation invariant state on $\clb=\otimes_{\IZ}\!M_d(\!C)$ and $\omega_R$ be it's restriction on $\clb_R=\otimes_{\IN}\!M_d(\!C)$. Then $\omega_R$ is either 
a type-I or type-III factor state.    
\end{thm}

\vsp 
In a related paper [Mo5] we deal with a complete weak classification of translation dynamics on $\clb=\otimes_{\IZ} \!M_d(\!C)$ with pure states with Kolmogorov's property [Mo1]. 
It says that $(\clb,\theta,\omega)$ and $(\clb,\theta,\omega')$ are unitarily equivalent 
for any two translation invariant states $\omega$ and $\omega'$ on $\clb$ if they satisfy Kolmogorov's property. In the proof simplicity of $\clb_R$ plays an important role and one can 
generalize replacing $\clb_R$ by a simple infinite dimensional separable $C^*$-algebra. [Mo4] address an asymptotic criteria for $\omega_R$ to be a type-I state.

\section{ Purity and Kolmogorov's property of stationary Markov processes: }

\bigskip
Let $\IT$ be either $\IR$, set of real numbers or $\IZ$ , set of integers and $\IT_+$ is non-negative numbers of $\IT$. Let $\tau_t:\clm_0 \raro \clm_0,\;\IT_+$ be a semi-group of completely 
positive unital normal maps on a von-Neumann algebra $\clm_0 \subseteq \clb(\clh_0)$ with a normal invariant state $\phi_0$. In case we are dealing with $\IT=\IR$, we also assume that for each 
$x \in \clm_0$, the map $t \raro \tau_t(x)$ is continuous in the weak$^*$ topology. Here we review the construction of stationary Markov process given as in [AM] in order to fix the notations 
and important properties.

\vsp
We consider the class $\ul{\clm}$ of $\clm_0$ valued functions
$\ul{x}: \IT \raro \clm_0$ so that $x_r \neq I$ for finitely
many supported points and equip with the point-wise multiplication
$(\ul{x}\ul{y})_r=x_ry_r$. We define the map $L: (\ul{\clm} \times \ul{\clm})
\raro \IC $ by
\be
L(\ul{x},\ul{y}) =
\phi_0(x_{r_n}^*\tau_{r_n-r_{n-1}}(x_{r_{n-1}}^*(.....x_{r_2}^*
\tau_{r_2-r_1}(x_{r_1}^*y_{r_1})y_{r_2})...y_{r_{n-1}})y_{r_n})
\ee
where $\ul{r}=(r_1,r_2,..r_n)\;r_1 \le r_2 \le .. \le r_n$ is
a finite collection of points in $\IZ$ containing both the support sets of 
$\ul{x}$ or $\ul{y}$. That this kernel is well defined
follows from our hypothesis that $\tau_t(I)=I, \; t \ge 0$ and
the invariance of the state $\phi_0$ for $\tau.$ The
complete positiveness of $\tau$ implies that the map $L$ is a
non-negative definite form on the set $\ul{\clm}$. Thus there exists a
Hilbert space $\clh$ and a map $\lambda: \ul{\clm} \raro \clh$ such
that $$<\lambda(\ul{x}),\lambda(\ul{y}) >= L(\ul{x},\ul{y}).$$
Often we will omit the symbol $\lambda$ to simplify our
notations unless more then one such maps are involved. For detailed and historical 
account we refer to articles [AFL,EL,Vi,Sa,BhP]. 

\vsp
We use the symbol $\Omega$ for the unique element in $\clh$
associated with $x=(x_r=I,\;r \in \IT )$ and the associated
vector state $\phi$ on $B(\clh)$ defined by $\phi(X)=<\Omega,X
\Omega>$.

\vsp
For each $t \in \IT$ we define shift operator $S_t: \clh \raro
\clh$ 
by the following prescription:
\be
(S_t\ul{x})_r = x_{r+t}
\ee
It is simple to note that $S = (( S_t ,\;t \in \IT))$ is a unitary
group of operators on $\clh$ with $\Omega$ as an invariant element and 
the map $t \raro S_t$ is continuous in strong operator topology as the
map $(t,x) \raro \tau_t(x)$ is sequentially continuous in weak$^*$ 
topology. For details we refer to [AM]. 

\vsp
For any $t \in \IT$ we set
$$\ul{\clm}_{t]}= \{\ul{x} \in \ul{\clm},\; x_r=I\;\forall r > t \} $$
and $F_{t]}$ for the projection onto $\clh_{t]}$, the closed
linear span of $\{\lambda(\ul{\clm}_{t]})\}$. For any $x \in \clm_0$
and $t \in \IR$ we also set elements $i_t(x),\in \ul{\clm}$ defined
by $$i_t(x)_r= \left \{ \begin{array}{ll} x ,&\; \mbox{if}\; r=t
\\ I,&\; \mbox{otherwise}\;  \end{array} \right.$$ 

\vsp
We also note that $i_t(x) \in \ul{\clm}_{t]}$ and set $\star$-homomorphisms
$j^0_0: \clm_0 \raro \clb(\clh_{0]})$ defined by
$$
j^0_0(x)\ul{y}= i_0(x)\ul{y}
$$
for all $\ul{y} \in \ul{\clm}_{0]}.$ That $j^0_0(x)$ is a well defined isometry follows
from (2.1) once we verify that it preserves the non-negative kernel whenever $x$ is 
an isometry. For any arbitrary element $x$ which is a linear sum of at most four isometry, we 
extend by linearity [AM]. Now we define $j^f_0: \clm_0 \raro \clb(\clh)$
by
$$
j^f_0(x)=j_0^0(x)F_{0]}.
$$
Thus $j^f_0$ is a $*$-homomorphism of $\clm_0$ at time $t=0$ with
$j^f_0(I)=F_{0]}$. Now we use the shift $(S_t)$ to obtain the
process $j^f=(j^f_t: \clm_0 \raro \clb(\clh),\;t \in \IT )$ and
forward filtration $F=(F_{t]},\;t \in \IT)$ defined by the
following prescription:
\be
j^f_t(x)=S_tj^f_0(x)S^*_t\;\;\;F_{t]}=S_tF_{0]}S^*_t,\;\;t \in \IT.
\ee

\vsp
So it follows by our construction that $ j^f_{r_1}(y_1)j^f_{r_2}(y_2)...
j^f_{r_n}(y_n) \Omega = \ul{y}$ where $y_r=y_{r_i},\;$ if $r=r_i$
otherwise $I,\;(r_1 \le r_2 \le .. \le r_n)$.  Thus $\Omega$ 
is a cyclic vector for the von-Neumann algebra $\clm_{[-\infty}$ generated 
by $\{ j^f_r(x), \;r \in \IT, x \in \clm_0 \}$. From (2.4) we also conclude that 
$S_tXS^*_t \in \clm_{[-\infty}$ whenever $X \in \clm_{[-\infty}$ and thus we can set 
a family of automorphism $(\al_t)$ on $\clm$ defined by $$\al_t(X)=S_tXS^*_t.$$ Since 
$\Omega$ is an invariant element for $(S_t)$, $\phi$ is an invariant normal state for $(\alpha_t)$.

\vsp
We also claim that
\be
F_{s]}j^f_t(x)F_{s]}=j^f_s(\tau_{t-s}(x))\;\; \forall s \le t.
\ee
For that purpose we choose any two elements $\ul{y},\ul{y'} \in
\lambda({\ul{\clm}}_{s]})$ and check the following steps with the aid of
(2.1) and (2.2): $$ <\ul{y},F_{s]}j^f_t(x)F_{s]}\ul{y'}>=
<\ul{y},i_t(x)\ul{y'}>$$
$$=<\ul{y},i_s(\tau_{t-s}(x))\ul{y'})>.$$ Since $\lambda
(\ul{\clm}_{s]})$ spans $\clh_{s]}$ it completes the proof of our claim. 

\vsp 
We set as in [Mo1] a tower of increasing von-Neumann algebras $\clm_{[t}$ with 
decreasing parameter $t \in \IT$ by setting $\clm_{[t}$ to be the von-Neumann algebra 
generated by $\{j^f_s(x): t \le s < \infty,\;x \in \clm_0 \}$ and so 
$\clm_{[-\infty}= \vee_{t \in \IT} \clm_{[t}$. Forward Markov property also 
ensures that $F_{t]}\clm_{[t}F_{t]}=j_t(\clm_0)$ [Mo1].

\vsp
For any element $\ul{x} \in \ul{\clm}$, we verify by the relation
$<\ul{y},F_{t]}\ul{x}=<\ul{y},\ul{x}>$ for all $\ul{y} \in {\cal M}_{t]}$
that $$(F_{t]}\ul{x})_r=
\left \{ \begin{array}{lll} x_r,&\;\mbox{if}\;r < t;\\
\tau_{r_k-t}(...\tau_{r_{n-1}-r_{n-2}}(\tau_{r_n-r_{n-1}}(
x_{r_n})x_{r_{n-1}})...x_t),&\;\mbox{if}\;r=t\\
I,&\;\mbox{if}\;r > t \end{array} \right. $$
where $r_1 \le ..\le r_k \le t \le .. \le r_n$ is the support of
$\ul{x}$. 
The result follows once we note that for any fix $\ul{x},\ul{y} \in \clh$
if $t \le r_1,r'_1$, where
$r_1,r'_1$ are the lowest support of $\ul{x}$ and $\ul{y}$ respectively,
$$<\ul{x},F_{t]}\ul{y}> = <F_{t]}\ul{x},F_{t]}\ul{y}>$$
$$=
\phi_0[ (\tau_{r_1-t}(...\tau_{r_{n-1}-r_{n-2}}(\tau_
{r_n-r_{n-1}}(x_{r_n})x_{r_{n-1}})...x_{r_1}) )^*
\;$$
$$\tau_{r'_1-t}(...\tau_{r'_{m-1}-r'_{m-2}}(\tau_
{r'_m-r'_{m-1}}(y_{r'_m})y_{r'_{m-1}})...y_{r'_1})]$$
Thus $F_{t]} \raro |\Omega><\Omega|$ as $t \raro -\infty$ if and only if 
$\phi_0(\tau_t(x)\tau_t(y)) \raro \phi_0(x)\phi_0(y)$ as $ t \raro
\infty$ for all $x,y \in \clm_0$. In such a case we have $\clm_{[-\infty}=\clb(\clh)$. Further since $S_t\Omega=\Omega$ and $F_{t]}\Omega=\Omega$, 
we get an imprimitivity systems $(\clh \ominus |\Omega><\Omega|, S_t, F_{t]} \ominus |\Omega><\Omega|:t \in \IT)$ for the group $\IT$, by 
restricting actions of $S_t,F_{t]}$ on the orthogonal complement $I_{\clh} \ominus |\Omega><\Omega|$ of $|\Omega><\Omega|$. Such a property is called 
[Mo1] {\it Kolmogorov's property } for $(\tau_t)$ and associated imprimitivity system $(I_{\clh} \ominus |\Omega><\Omega|,S_t,F_{t]} \ominus |\Omega><\Omega|)$ 
is called {\it Kolmogorov shift }. We will deal with this notion in the context of quantum spin chain in section 4 and 5. However Kolmogorov's property is not necessary 
for $\clm_{[-\infty}$ to be $\clb(\clh)$ even when $\phi_0$ is faithful [Mo**]. We will not use this in this paper. 

\vsp
We also recall that triplet $(\clm_0,\tau_t,\phi_0)$ is called ergodic, strong mixing if $(S_t)$ is ergodic, strong mixing respectively. For details we refer to [AM] 
where simple correspondence were established. For further details we also refer to [Mo1]. 
 
\vsp 
The following theorem gives a necessary and sufficient condition for purity i.e. $\clm_{[-\infty}=\clb(\clh)$.  

\vsp 
\begin{thm} The following statements are equivalent: 

\NI (a) For each $t \ge 0$, there exists an element $x_t \in \clm_0$ so that $x_t \raro I$ in strong 
operator topology and for each $s \ge 0$, $\phi_0(\tau_t(x)x^*_{s+t}x_{s+t}\tau_t(y)) \raro \phi_0(x)\phi_0(y)$ 
as $t \raro \infty$ for all $x,y \in \clm_0$;

\NI (b) $\vee_{t \in \IR } \clm_{[t} = \clb(\clh)$ 

\end{thm} 

\begin{proof} 
First we will prove (a) implies (b). Fix such a sequence $x_t$ with property (a). We consider the sequence 
of elements $j_{-t}(x_t): t \ge 0 \}$. We claim that $j_{-t}(x_t) \raro |\Omega><\Omega|$ as $t \raro \infty$ in strong operator 
topology. To that end we fix an element $\ul{y}$ with support 
$(s_0 \le s_1 \le ..\le s_n)$ and $y(s_k)=y_k$ to check that for all $-t \le s_0$ we have 
$$||j_{-t}(x_t)\ul{y}||^2 = \phi_0(\tau_{t+s_0}(y^*) x^*_tx_t\tau_{t+s_0}(y))$$ 
where $y=\tau_{s_n-s_{n-1}}( y_n\tau_{s_{n-1}-s_{n-2}}( y_{n-2} \tau_{s_1-s_0}( y_1) y_0))$
as 
$$F_{s_0}j_{s_n}(y_n)j_{s_{n-1}}(y_{n-1})..j_{s_0}(y_0)F_{s_0}$$
$$=j_{s_0}( \tau_{s_n-s_{n-1}}( y_n\tau_{s_{n-1}-s_{n-2}}( y_{n-2} \tau_{s_1-s_0}( y_1) y_0) ) )$$
Thus $$||j_{-t}(x_t)\ul{y}||^2 = $$
$$ \phi_0(\tau_{t+s_0}(y^*)x_t^*x_t\tau_{t+s_0}(y))$$ 
Thus    
$$ \raro |\phi_0(y)|^2 = |<\Omega,\ul{y}>|^2$$ as $t \raro \infty$ by (a) for all $s_0 \le 0$.    

Similarly we also have 
$$<\ul{y'},j_{-t}(x_t)\ul{y}> = <j_{-t}(I)\ul{y'},j_{-t}(x_t)\ul{y}> $$
$$= \phi_0(\tau_{t+s_0}(y'^*)x_t \tau_{t+s_0}(y)) \raro 
<\ul{y'},\Omega><\Omega,\ul{y}>$$ as $t \raro \infty$ follows along the same line since for each $s \ge 0$,
$$\phi_0(\tau_t(x)x_{s+t}\tau_t(y)) \raro \phi_0(x)\phi_0(y)$$ for all $x,y \in \clm_0$ as $t \raro \infty$ ( as first we 
use Cauchy-Schwarz inequality to prove the statement for all $x$ for which $\phi_0(x)=0$ and then replace $x$ by $x-\phi_0(x)I$ ).    
We combine now above two statements to prove $j_{-t}(x_t) \raro |\Omega><\Omega|$ in strong operator topology as $t \raro \infty$ 
as the family $\{j_{-t}(x_t): t \ge 0 \}$ is uniformly norm bounded as $x_t \raro 1$ in strong operator topology by (a) and each $j_{-t}$ 
is an injective isomorphism (in particular contractive property ).     

\vsp 
For the converse we use Kaplansky's density theorem to ensure existence of contractive elements $Y_t \in \clm_{[t}$ so that 
$Y_t \raro |\Omega><\Omega|$ in strong operator topology as $t \raro -\infty$. Since 
$F_{t]}\clm_{[t}F_{t]} = j_t(\clm_0)$ we also have $F_{t]}Y_{[t}F_{t]} = j_t(y_t)$ for some $y_t \in \clm_0$. 
Since $F_{t]} \raro F_{-\infty]} \ge |\Omega><\Omega|$ in strong operator topology as $t \raro -\infty$, we 
conclude also that $j_t(y_t) \raro |\Omega><\Omega|$ in strong operator topology. 
Now we compute for $t \ge s \ge 0$ that 
$$<j_{-t}(y_{-t})j_{-s}(y)\Omega,j_{-t}(y_{-t})j_{-s}(x)\Omega>=\phi_0(\tau_{t-s}(y^*)y_{-t}^*y_{-t}\tau_{t-s}(x))$$ and thus we conclude (a) with $x_t=y_{-t}$ for all $t \ge 0$. 
Taking $x=y=1$, we get $x_t\Omega \raro \Omega$ as $t \raro \infty$ strongly and thus by separating property of $\Omega$, we get $x_tf \raro f$ strongly for a dense set and thus
$x_t \raro I$ in strong operator topology being a uniformly bounded family.  \end{proof} 

\vsp 
That Kolmogorov's property implies strong mixing [AM,Mo1] i.e. $\tau_t(x) \raro \phi_0(x)$ in weak operator topology 
follows by Cauchy-Schwarz inequality 
$$|\phi_0(\tau_t(x^*)y)| \le \phi_0(\tau_t(x^*)\tau_t(x))\phi_0(y^*y)$$ 
for $x,y \in \clm_0$. Within the frame work of classical probability theory it is well known that Kolmogorov's property is preserved when we reverse the direction of time evolution of an automorphism on a measure space preserving a probability measure. A simple proof follows from Kolmogorov-Sinai-Rohklin theorem which gives criteria of positive dynamical entropy [Pa] for an automorphism to be a Kolmogorov's automorphism. Such a notion in terms of a quantum dynamical entropy is still missing within the frame work of Connes-St\o rmer dynamical entropy [NeS]. Further now one can as well ask the same question about weak Kolmogorov's or purity property i.e. whether such a property is time reversible. In the following text we aim to investigate how time reverse process is related to the notion of weak Kolmogorov's 
property and also with Kolmogorov's property. Following [AM], we will consider the time reverse process associated with the KMS-adjoint quantum dynamical semi-group $(\clm_0,\tilde{\tau},\phi_0)$. 
First we recall from [AM] time reverse process associated with the KMS-adjoint semi-group in the following paragraph.

\vsp
Let $\phi_0$ be also a faithful state and without loss of generality let also 
$(\clm_0,\phi_0)$ be in the standard form $(\clh_0,\clm_0,\clj,\Delta,{\cal P},\omega_0)$
[BR] where $\omega_0 \in \clh_0$, a cyclic and separating vector for $\clm_0$, so
that $\phi_0(x)= <\omega_0,x\omega_0>$ and the closure of the close-able operator 
$S_0:x\omega_0 \raro x^*\omega_0, S$ possesses a polar decomposition
$S=\clj \Delta^{1/2}$ with the self-dual positive cone $\clp$ as the closure 
of $\{ \clj x \clj x\omega_0:x \in \clm_0 \}$ in $\clh_0$. Tomita's [BR] theorem says 
that $\Delta^{it}\clm_0\Delta^{-it}=\clm_0,\;t
\in \IR$ and $\clj \clm_0 \clj=\clm'_0$, where $\clm'_0$ is the
commutant of $\clm_0$. We define the modular automorphism group
$\sigma=(\sigma_t,\;t \in \IT )$ on $\clm_0$
by
$$\sigma_t(x)=\Delta^{it}x\Delta^{-it}$$ which satisfies the modular  relation
$$\phi_0(x\sigma_{-{i \over 2}}(y))=\phi_0(\sigma_{{i \over 2}}(y)x)$$
for any two analytic elements $x,y$ for the automorphism. A more useful 
form for modular  relation here
$$\phi_0(\sigma_{-{i \over 2}}(x^*)^* \sigma_{-{i \over 2}}(y^*))=\phi_0(y^*x)$$ 
which shows that $\clj x\Omega= \sigma_{-{i \over 2}}(x^*)\Omega$. Furthermore for 
any normal state $\psi$ on $\clm_0$ there exists an unique vector $\zeta \in {\cal P}$ 
so that $\psi(x)= <\zeta,x\zeta>$.

\vsp
We consider the unique Markov semi-group $(\tau'_t)$ on the commutant
$\clm'_0$ of $\clm_0$ so that $\phi_0(\tau_t(x)y)=\phi_0(x\tau'_t(y))$ for all
$x \in \clm_0$ and $y \in \clm'_0$. Proof follows a standard application of Dixmier 
lemma a variation of Radon-Nikodym theorem [OP]. We define weak$^*$ continuous
Markov semi-group $(\tilde{\tau}_t)$ on $\clm_0$ by $\tilde{\tau}_t(x)=\clj \tau'_t(\clj x \clj ) \clj.$
Thus we have the following adjoint relation
\be
\phi_0(\tau_t(x)\sigma_{-{i \over 2}}(y))=\phi_0(\sigma_{i \over 2}(x) \tilde{\tau}_t(y))
\ee
for all $x,y \in \clm_0$, analytic elements for $(\sigma_t)$. 

\vsp
We consider the forward weak Markov processes $(\clh,S_t,j^f_t,F_{t]},F_{[t},\;\;t \in \IT,\;\Omega)$ 
associated with $(\clm_0,\tau_t,\;t \ge 0,\;\phi_0)$ and the forward weak Markov processes 
$(\tilde{H},\tilde{S}_t,\;\tilde{j}^f_t,\tilde{F}_{t]},\tilde{F}_{[t},,\; t \in \IT,\; 
\tilde{\Omega})$ associated with $(\clm_0,\tilde{\tau}_t,\; t \ge 0,\;\phi_0)$. 

\vsp
Now following a basic idea of G.H. Hunt, as in [AM], here we consider the transformation 
$\clh \raro \tilde{\clh}$ generalizing Tomita's conjugate operator given by 
$\tilde{x}(t) = \sigma_{-{i \over 2}}(x(-t)^*)$ for $x \in \clh$ supported on analytic elements of the modular 
automorphism group of $\phi_0$. We recall that modular condition and duality property will ensure that such a 
transformation is anti-inner product preserving and thus extends to an anti-unitary operator 
$U_0:\clh \raro \tilde{\clh}$ which takes $x$ to $\tilde{x}$. Further there exists a unique backward weak Markov 
processes $(j^b_t),(\tilde{j}^b_t)$ which generalizes Tomita's representation $\pi^b_0:x \raro \clj_0 \pi_0(x^*) \clj_0 \in \clm_0'$ ( linear anti-isomorphism ) so that 
$$F_{[s}j^b_t(x)F_{[s}=j^b_s(\tilde{\tau}_{s-t}(x))$$ 
for $-\infty < t \le s < \infty $. 
$$\tilde{F}_{[s}\tilde{j}^b_t(x)\tilde{F}_{[s}=\tilde{j}^b_s(\tau_{s-t}(x))$$ 
for $-\infty < t \le s < \infty $. We set tower of increasing von-Neumann algebras $\clm^b_{t]}$ in backward 
direction with increasing parameter $t \in \IT$ by setting 
$\clm^b_{t]}= \{j^b_s(x): -\infty < s \le t, \; x \in \clm_0 \}''$ 
and $\clm^b_{\infty]} = \vee_{ t \in \IT } \clm^b_{t]}$. Once more we have 
$F_{[t}\clm^b_{t]}F_{[t}=j^b_t(\clm_0)$ for each $t \in \IT$ by backward Markov property [Mo1].   

\vsp 
We have more details in the following theorem.  

\vsp
\begin{pro} 

[AM]  There exists an unique anti-unitary operator $U_0:\clh \raro \tilde{\clh}$ so that 

\NI (a) $U_0 \Omega = \tilde{\Omega}$;

\NI (b) $U_0 S_t U^*_0 = \tilde{S}_{-t}$ for all $t \in \IT$;

\NI (c) $U_0F_{t]}U^*_0=\tilde{F}_{[-t},\;\;U_0F_{[t}U^*_0=\tilde{F}_{-t]}$ 
for all $t \in \IT$;

\NI (d) We set backward processes $j^b_t:\clm_0 \raro \clb(\clh)$ and $\tilde{j}^b_t:\clm_0 \raro \clb(\tilde{\clh})$
as injective $*$-anti-homomorphism so that  
$$j^b_t(x) = U^*_0 \tilde{j}^f_{-t}(x^*) U_0$$ 
$$\tilde{j}^b_t(x)= U_0j^f_{-t}(x^*) U^*_0 $$ 
for all $t \in \IT$; Then we have backward Markov property for $s \le t$:
$$F_{[t}j^b_s(x)F_{[t}=j^b_t(\tilde{\tau}_{t-s}(x)),\;\;,
\tilde{F}_{[t}\tilde{j}^b_s(x)\tilde{F}_{[t}=\tilde{j}^b_t(\tau_{t-s}(x))$$
Further $F_{[t}\clm^b_{t]}F_{[t}=j^b_t(\clm_0)$ and $\tilde{F}_{[t}\tilde{\clm}^b_{t]}\tilde{F}_{[t}=\tilde{j}^b_t(\clm_0)$ for all $t \in \IT$.   

\end{pro} 

\begin{proof} 
Properties (a) -(c) are evident by our construction of anti-unitary operator $U_0$. Backward 
weak Markov property follows once we use $U_0$ for push forward method from forward weak Markov property of 
dual forward Weak Markov processes. Last part also follows via $U_0$ as we have 
$F_{t]}\clm_{[t}F_{t]}=j^f_t(\clm_0)$ and $\tilde{\tilde{\tau}}_t=\tau_t;t \ge 0$. \end{proof}

\vsp 
Exploring Markov property of both forward and backward processes, Tomita's duality also ensures the following 
Hunt duality theorem.  

\vsp
\begin{pro} [Mo2] For each $t \in \IZ$, $\clm'_{[t}=\clm^b_{t]}$, where $\clm_{[t}=\{j^f_s(x):\;x \in \clm_0, 
s \ge t \}''$ and $\clm^b_{t]}=\{j^b_s(x): x \in \clm_0,\; s \le t \}''$. 
\end{pro}

\vsp 
\begin{proof} 
We refer to Proposition 3.7 in [Mo2]. Note that it only needs weak Markov property for both forward and
backward processes and Tomita's duality relation at fiber at $t \in \IT$. 
\end{proof}

\vsp 
\begin{thm}
The following statements are equivalent: 

\NI (a) For each $t \ge 0$, there exist a contractive element $x_t \in \clm_0$ so that for each $s \ge 0$, $x_{s+t}\tau_t(x) \raro \phi_0(x)$ in strong operator topology i.e. 
$\phi_0(\tau_t(x)x^*_{s+t}x_{s+t}\tau_t(y)) \raro \phi_0(x)\phi_0(y)$ as $t \raro \infty$ for all $x,y \in \clm_0$;

\NI (b) $\vee \clm_{[t} = \clb(\clh)$;

\NI (c)  $||\psi \tilde{\tau}_t-\phi_0|| \raro 0$ as $t \raro \infty$ for any normal state $\psi$ on $\clm_0$;

\NI (d) $\bigcap \clm^b_{t]}= \IC$ 
\end{thm} 

\vsp 
\begin{proof} 
 That (a) and (b) are equivalent follows from Theorem 2.1 as separating property of $\Omega$ for 
$\clm_0$ will ensure that $\Omega$ is cyclic for commutant of $\clm_0$ and thus $x_{s+t}\tau_t(x)\Omega \raro \phi_0(x)\Omega$ 
in strong topology is equivalent to $x_{s+t}\tau_t(x) \raro \phi_0(x)1$ in strong operator topology. That (c) and (d) are equivalent 
follows by Proposition 1.1 in [Ar1] and Theorem 2.4 in [Mo2] as support projection of the state in $\clm^b_{0]}$ 
is $[(\clm^b_{0]})'\Omega]=[\clm_{[0}\Omega]=F_{[0}$ and using time-reversal operator $U_0$ we also check that for 
any $Y \in \clm^b_{0]}$ there exists a $y \in \clm_0$ so that $j^b_0(\tilde{\tau}_t(y)) = F_{[0}\alpha_{-t}(Y)F_{[0}$. 
That (b) and (d) are equivalent follows from duality proved in Proposition 2.3. 
\end{proof}
\begin{thm}  
Let $\clm_0$ be also a factor as in Theorem 2.4. If $\psi_{[-\infty}$ is pure then $\clm_0$ is not a type-II factor.  
\end{thm}  
\begin{proof}  
We will rule out the possibility for $\clm_0$ to be of type-II factor by bringing a contradiction. Let $\clm_0$ be type-II and $tr_0$ be a 
semi-finite weight in case $\clm_0$ is type-II${_\infty}$ otherwise unique tracial state. In any case for each $t \in \IT$, $\clm_{[t}$ is also 
a type-II factor by Proposition 4.1 in [Mo3]. $j_0$ being faithful, $j_0(p)$ is also a finite projection in $\clm_{[-t}$ for each $t \ge 0$ for 
a finite non-zero projection $p$ in $\clm_0$. We consider the nested family 
of increasing type-$II_1$ factors $\{\clm^p_t=j_0(p)\clm_{[-t}j_0(p): t \ge 0 \}$ 
and also set $\clm_t=j_0(I)\clm_{[-t}j_0(I)$. By the uniqueness of a normal tracial state on a finite von-Neumann algebra [Di,chapter 4] we get a tracial functional $tr^p_0$ 
on $*$-sub-algebra $\clm^0(p) = \bigcup_{t \ge 0} \clm^p_t$ such that $tr^p_0(p)=tr_0(p)$. Thus $tr^p_0$ has an unique extension on the $\clm^1(p)$, 
$C^*$-norm closure of $\clm^0(p)$ so that $tr^p_0(p)=tr_0(p)$. $tr^p_0$ is also a trace on $\clm^1(p)$. 

\vsp 
We set notation $\clm^0=\bigcup_{t \le 0} \clm_t$, $\clm^1$ for it's $C^*$-completion and $\clm$ for it's von-Neumann completion i.e. $\clm= \vee_{t \le 0} \clm_t$. 
The family of traces $\{tr^p_0: tr_0(p) < \infty \}$ induces a tracial weight $tr_0$ on $\clm^1$. We consider GNS representations $\pi_0:\clm^1 \raro \clb(\clh_{\pi_0})$ 
of the semi-finite weight $tr_0$ on $\clm^1$. We consider also GNS space $(\clh^p_{\pi_0},\pi^p_0,\zeta_{tr^p_0})$ associated with $tr^p_0$ on $\clm^1(p)$. For two 
finite projection $p \le p'$ in $\clm_0$, GNS space $\clh^p_{\pi_0}$ has a natural embedding into $\clh^{p'}_{\pi_0}$ as a closed subspace and the Hilbert space $\clh_{\pi_0}$ can as well 
be realized as the inductive limit Hilbert space of the family $\{ \clh^p_{\pi_0}:p,\;tr_0(p) < \infty \}$ with respect to the natural inclusion 
map and tracial weight $tr_0$. We claim that $tr^p_0$ has a unique normal extension to $(\clm^1(p))''=\clm(p)$, where $\clm(p)=\vee_{t \le 0} \clm^p_t$. 

\vsp 
For $t \le 0$ we set $k_t(x)= \pi_0(j_t(x)) $ for all $x \in \clm_0$ and $t \le 0$ for which $tr_0(x^*x) < \infty$ and 
set $\cln^p_t=\pi_0(\clm^p_t)$ for $t \le 0$. Note that $k_t(I)$ which may not be defined via homomorphism. We defined $k_t(I)$ to be the 
limit of the net $k_t(p)$ of increasing projections of finite projections $p \uparrow I$ in $\clm_0$. It is clear by $*$-homo-morphism property 
of $\pi_0$ that $x \raro k_t(x)$ is also (weak) Markov process associated with $(\clm_0,\tau_t)$ i.e. 
$$k_s(I)k_t(x)k_s(I)=k_{s}(\tau_{t-s}(x))$$ 
for all $s \le t \le 0$ and $x \in \clm_0,\;tr_0(x^*x) < \infty$. We choose a net $p_{\alpha} \uparrow I$ of finite
projections. By Cauchy-Schwarz inequality for any normal state $\phi$ on $\clm$ we have 
$$|\phi(k_t(p_{\alpha}\tau_{t-s}(x)(I-p_{\alpha}))| \le ||x|| \phi(k_t(I-p_{\alpha})) \raro 0$$ 
as $p_{\alpha} \uparrow I$. Further $\clh_{\pi_0}$ being the GNS space, $\zeta_{tr_0}$ is also cyclic for 
$\{\pi_0(X):X \in \clm^0 \}$ and thus $\zeta_{tr_0}$ is also cyclic for $\{k_t(x): t \le 0, tr_0(x^*x) < \infty \}$. 

\vsp 
First we consider the simplest situation namely $\clm_0$ is a type-II$_1$ factor and $p=1$ and $\phi_0$ is the unique trace $tr_0$ on $\clm_0$. We identify $\clm_0$ with it's 
standard form associated with the trace $tr_0(x)=<\zeta_{tr_0},x \zeta_{tr_0}>$. In such a case we set unitary operator $U:\clh_{0]} \raro \clh_{\pi_0}$ given by 
\be 
U: j_{t_1}(x_1)j_{t_2}(x_2)...j_{t_n}(x_n) \Omega \raro k_{t_1}(x_1)k_{t_2}(x_2)...k_{t_n}(x_n) \zeta_{tr_0} 
\ee 
for $t_n \le t_{n-1} ..\le t_1 \le 0$ and any elements $x_k \in L^2(\clm_0,tr_0)$.  

More generally for a faithful normal state $\phi_0$ on $\clm_0$ we have an unique vector $\zeta_{\phi_0}$ in the self dual cone associated with the normalize trace 
such that $\phi_0(x)=<\zeta_{\phi_0},x\zeta_{\phi_0}>$. In such a case we modify $U$ by replacing $\zeta_{tr_0}$ by $\zeta_{\phi_0}$, where we note that the GNS space associated 
with the tracial state is a closed subspace of the GNS space $\clh_{\pi_0}$ and the modification is well defined. That $U$ is indeed an inner product preserving map follows using 
Markov property (2.4). So by our construction we have $Uj_t(x)U^*=k_t(x)$ for all $x \in \clm_0$ and $t \le 0$. Thus von-Neumann algebra 
$\clm = \vee_{t \le 0} \clm_t$ is isomorphic to von-Neumann algebra $\cln = \vee_{t \le 0} \cln_t$. 

\vsp 
For a more general situation with type-II$_{\infty}$ factor, we choose
\be 
U_p: j_0(p)j_{t_1}(x_1)j_{t_2}(x_2)...j_{t_n}(x_n)j_0(p)\Omega \raro k_0(p)k_{t_1}(x_1)k_{t_2}(x_2)...k_{t_n}(x_n)k_0(p)\zeta_{tr_0} 
\ee 
for $t_n \le t_{n-1} ..\le t_1 \le 0$ and any elements $x_k \in L^2(\clm_0,tr_0)$ where $\phi_0(x)=<\zeta_0,x\zeta_0>$ where $\zeta_0$ 
is the unique unit vector in the positive cone associated with the semi-finite weight $tr_0$ [Ha].  

\vsp 
Thus we may assume without loss of generality that $\clm(p)=\pi_0(\clm^1(p))''$ and $\clh_{0]}$ is the GNS space associated with $tr_0$ on $\clm^1(p)$. Thus $tr_0$ has a unique normal extension 
from $\pi_0(\clm_1(p))$ to $\clm(p)$ given by the vector state $\zeta_{tr_0}$. The vector state being normal, 
we conclude that $tr_0$ given by the vector state $\zeta_{tr_0}$ is also a trace on $\clm$.                   

\vsp 
However by Markov property we have $j_0(p) (\vee_{t \ge 0} \clm_{[-t} ) j_0(p) = \clm(p) $. Thus by purity 
$\clm(p)$ is equal to the set of all bounded operators on $j_0(p)$. Thus $j_0(p)$ and so $p$ is a finite dimensional 
subspace. This contradicts our starting assumption that $p \neq 0$ and is a finite projection in a type $II$ factor. 

\end{proof}  
 
\begin{thm}  
Let $(\clm_0,\tau_0,\phi_0)$ be as in Theorem 2.4 and $\vee \clm_{[t}=\clb(\clh)$. Then $\clg=^{\mbox{def}}\{x: \tilde{\tau}_t\tau_t(x)=x:\;t \ge 0 \}$ which is equal to $\{x: \tau_t(x^*)\tau_t(x)=\tau_t(x^*x),\tau_t(x)\tau(x^*)=\tau_t(xx^*),\;\tau_t \sigma_s(x)
=\sigma_s \tau_t(x):\;t \ge 0, s \in \IR \}$ is trivial. 
\end{thm} 

\begin{proof} 
We consider the dynamics $(\clg,\tau_t,\phi_0)$. By Proposition 2.2 (a) in [Mo2] there exists a conditional expectation $E:\clm \raro \clg$ onto $\clg$ and modular group $(\sigma_t:t \in \IR)$ commutes with 
$(\tau_t:t \ge 0)$ on $\clg$. We claim that $(\clg_0,\tau_t,\phi_0)$ satisfies Kolmogorov's property. As a first 
step we verify that $\phi_0(\clj y \clj \tilde{\tau}_t \tau_t(x)) = \phi_0(\clj \tau_t(y) \clj \tau_t(x)) \raro \phi_0(y)\phi_0(x)$ for all $x,y \in \clm_0$ as follows: 

We fix any $y \ge 0$ so that $\phi_0(y)=1$ and consider the normal state $\phi(x)=\phi_0(\clj y \clj x)$ on 
$\clm_0$ and compute the following simple steps 
$$|\phi \tilde{\tau}_t \tau_t(x)-\phi_0(x)| \le ||\phi \tilde{\tau}_t -\phi_0|| ||\tau_t(x)|| $$
$$\le ||\phi \tau_t-\phi_0||\;||x|| \raro 0$$ 
as $t \raro \infty$. For more general $y$ we write it as a linear combination of possibly four such 
non-negative elements and use linearity of the maps involved to complete the proof of the claim. We write now
$\phi_0(\clj \tau_t(x) \clj \tau_t(x))=||\Delta^{1/4}\tau_t(x^*)\Omega||^2$ and note that we have shown that
$\Delta^{1 \over 4}\tau_t(x^*)\Omega \raro \phi_0(x^*) \Omega$ strongly as $t \raro \infty$. So far $x \in \clm$ but
now we will restrict $x$ to $\clg$. Since modular group $(\sigma_s:s \in \IR)$ commute with $(\tau_t:t \ge 0)$ on 
$\clg$, weak$^*$ dense subspace of analytic elements [BR1, Proposition 2.5.22 ] of the form $\{ x_n= 
\sqrt{n \over \pi} \int \sigma_s(x)e^{-ns^2}ds: n \ge 1, x \in \clg \}$ are also preserved by $(\tau_t:t \ge 0)$. 

\vsp 
Thus we legitimately write taking $x^*$ instead $x$ 
$\Delta^{1 \over 4} \tau_t(x)\Omega=\Delta^{1 \over 4} \tau_t(x)\Delta^{-{1 \over 4}}\Omega=
\tau_t(\sigma_{i \over 4}(x))\Omega$ for such analytic elements $x \in \clg$ for $(\sigma_s)$. 
Thus we conclude that $\tau_t(x) \raro \phi_0(x)$ in strong operator topology for a weak$^*$ dense 
sub-space of $\clg$ of analytic elements $\clg_a$ of modular group $(\sigma_s)$. Now we remove the 
restriction on $x$ to be an analytic element as $[\clg_a\Omega]=[\clg\Omega]$ and $x\Omega \raro 
\tau_t(x)\Omega$ is a family of contraction on the GNS space associated with $\phi_0$.  

\vsp 
However since $\tau_t(x^*)\tau_t(x)=\tau_t(x^*x)$ for all $t \ge 0$ for $x \in \clg$, by invariance of the state we have 
$\phi_0(x^*x)=|\phi_0(x)|^2$ and thus faithful property of $\phi_0$ says that $x=\phi_0(x)I$. This completes the proof that $\clg$ is trivial 

\end{proof}  

\vsp 
Now we focus on the question: how situation changes when $\clm_0$ is a type-I factor or more generally type-I von-Neumann algebra with completely atomic center. 

\vsp 
\begin{thm}  
Let $(\clm_0,\tau_t,\phi_0)$ be as in Theorem 2.4 and $\clm_0$ be also a type-I von-Neumann algebra with center completely atomic. 
Then strong mixing property of $(\clm_0,\tau_t,\phi_0)$ is equivalent to all the statements (a)-(d) of Theorem 2.4.  
\end{thm} 

\begin{proof} 
Strong mixing property is time reversal i.e. $(\clm_0,\tilde{\tau}_t,\phi_0)$ is also strongly mixing. 
By a well known theorem of Dell'Antonio [De] $\clm_0$ being type-I with completely atomic center weak limit of a sequence 
of normal states in the pre-dual space of $\clm_0$ is equivalent to strong convergence thus strong mixing property 
$\psi \tilde{\tau}_t(x) \raro \phi_0(x)$ for each $x \in \clm_0$ implies in particular that $||\psi \tilde{\tau}_t - \phi_0|| 
\raro 0$ as $t \raro \infty$ for any normal state. \end{proof} 

\vsp 
In [Theorem 4.7 Mo1] a faulty proof appeared for a claim that strong mixing property of $(\clm_0,\tau_t, \phi_0)$ implies Kolmogorov's property for type-I von-Neumann algebra 
$\clm_0$ with completely atomic center. The argument used Proposition 4.6 in [Mo1] which is faulty and such a statement is not true unless modular group $(\sigma_s:s \in \IR)$ 
associated with $\phi_0$ commutes with $(\tau_t:t \in \IT_+)$. Though Proposition 4.6 in [Mo1] 
(if part of the statement) can not be repaired as the statement itself is faulty as argued above, in the following we give a correct proof for the statement 
given in Theorem 4.7 in [Mo1]. The proof that follows is far more involved then what one would expects at first sight.

\begin{thm} 
Let $\clm_0$ be also a type-I von-Neumann algebra acting on a separable Hilbert space $\clh_0$ with center completely atomic and 
$(\clm_0,\tau_t,\phi_0)$ is given as in Theorem 2.4. Then the following statements are equivalent:

\NI (a) $(\clm_0,\tau_t,\phi_0)$ is strong mixing property;

\NI (b) $(\clm_0,\tau_t,\phi_0)$ is Kolmogorov; 

\NI (c) $\{x: \tau_t(x^*)\tau_t(x)=\tau_t(x^*x),\tau_t(x)\tau(x^*)=\tau_t(xx^*):\;t \ge 0 \}=\{zI; \;z \in \IC \}$.  
\end{thm}  

\begin{proof}  
Once again we assume without loss of generality that $\clm_0$ is in standard form associated with $\phi_0$. Separating property of $\Omega$ for $\clm_0$ ensures that the strong mixing property of $(S_t)$ 
is equivalent to convergence of $\tau_t(x) \raro \phi_0(x)I$ as $t \raro \infty$ in weak operator topology for each $x \in \clm_0$ [AM]. Whereas Kolmogorov property is equivalent to convergence of $\tau_t(x) \raro \phi_0(x)I$ 
as $t \raro \infty$ in strong operator topology. We claim under our assumption on $\clm_0$ these two convergence are equivalent. We will make use of dilation described in Theorem 2.1. A direct proof is 
not visible at this moment which is desirable. As a first step of the proof, we first consider the case where $\IT$ is $\IZ$ i.e. time is discrete. 

\vsp 
We identify $\clm_0$ to direct sum $\sum_{ k \ge 1} \clb(\clh_k)$ where each $\clh_k$ is a complex separable Hilbert space. Let $p$ be a finite 
dimensional projection on $\clm_0$. $p$ being a compact operator, weak convergence of $\tau_n(x)$ ensures that $p\tau_n(x) \raro p \phi(x) $ as $n \raro \infty $ 
in strong operator topology. We claim that $j_{-n}(p) \raro \phi(p) |\Omega<\Omega|$ 
as $n \raro \infty$ in weak operator topology. The claim follows trivially once we recall by (2.4)    
$$<\lambda(\ul{x}),j_{-n}(x)\lambda(\ul{y})> $$
$$=<j_{m_k}(x_{m_k})j_{m_{k-1}}(x_{m_{k-1}})...j_{m_1}(x_{m_1})\Omega, j_{-n}(x)j_{m_k}(y_{m_k})j_{m_{k-1}}(y_{m_{k-1}})...j_{m_1}(y_{m_1})\Omega>$$
$$=<j_{-n}(I)j_{m_k}(x_{m_k})...j_{m_1}(x_{m_1})j_{-n}(I)\Omega,j_{-n}(x)j_{-n}(I)j_{m_k}(y_{m_k})...j_{m_1}(y_{m_1})j_{-n}(I)\Omega>$$
$$=<j_{-n}(x')\Omega, j_{-n}(x)j_{-n}(y')\Omega>$$
$$= \phi(\tau_{m_1+n}(x'^*)x\tau_{m_1+n}(y'))$$
where finite support of $\ul{x},\ul{y}$ is $(m_1,m_2,...,m_k)$ and $n \ge -m_1$ with elements  
$x'=x_{m_1}\tau_{m_2-m_1}(x_{m-2}\tau_{m_3-m_2}(.....))$ and $y'=y_{m_1}\tau_{m_2-m_1}(y_{m-2}\tau_{m_3-m_2}(.....))$ respectively. 

\vsp 
Let $\{p_k:k \ge 1 \}$ be a sequence of finite dimensional increasing projections in $\clh_1$ so that $p_k \raro I$ as $k \raro \infty$ in strong 
operator topology ( such a sequence exists as $\clm_0$ is type-I with center completely 
atomic ). We claim that there exists a sub-sequence $n_k:k \ge 1$ of natural numbers so that 
$$j_{-{n_k}}(p_k) \raro |\Omega><\Omega|$$ as $n \raro \infty$ in strong operator topology. 
It is enough if we show convergence in weak operator topology as limit is also a projection. 

\vsp 
Fix any vector $f \in \clh$ i.e. in the dilated space described in Theorem 2.1. For each $k$, $j_{-n}(p_k) \raro \phi(p_k)|\Omega<\Omega|$ in weak operator topology as $n \raro \infty$, we can subtract 
a sub-sequence $n_k(f):k \ge 1$ so that $<f,|\Omega><\Omega|f> - <f, j_{-n_k(f)}(p_k)f> \raro 0$ as $k \raro \infty$ as $p_k \uparrow I$ in strong operator topology by our choice. Now we can extract a 
sub-sequence $n_k$ independent of the vector $f$ for all $f$ in a countable dense subset of $\clh$ by Cantor's method such that 
$<f,|\Omega><\Omega|f> -<f,j_{-n_k}(p_k)f> \raro  0$ for all $f$ in a countable dense set for $\clh$, which is separable as $\clh_0$ is so (Proposition 2.7 (a) ). 
Since the family of operators involved are uniformly bounded, we conclude that $j_{-n_k}(p_k) \raro |\Omega><\Omega|$ in weak and so strong operator topology.      

\vsp 
We recall $F_{-\infty]} =s.\mbox{lim}_{n \raro \infty}F_{-n]}$ and $F_{-\infty]}  \ge |\Omega><\Omega|$. We claim now that $F_{-\infty]} =|\Omega><\Omega|$. Suppose not. Then we will have an unite vector $f \in \clh$ 
such that $<f,\Omega>=0$ and $|f><f| \le F_{-\infty]} $. We claim that there exists a sequence of positive contractive elements $y_n \in \clm_0=\clb(\clh_1)$ of finite rank such that 
$j_{-n}(y_n) \raro |\Omega><\Omega|$ as $n \raro \infty$ in strong operator topology. 

\vsp 
To that end we consider the $*$-sub-algebra $\clm^F_{[n}$ generated by elements 
$\{ j_m(x): m \ge n, x \in \clm^F_0 \}$ where $\clm_0^F$ is the $*$-algebra generated be elements of $\clm_0$ of finite rank. For any element $X \in \clm^F_{[n}$, we claim that  
$F_{n]}XF_{n]} \in j_n(\clm^F_0)$. For a proof we first note that $\tau$ being a linear map on $\clm_0$, it takes finite rank operator $x \in \clm_0$ to another finite rank 
operator and finite rank operators form an $*$ sub-algebra. Thus for an element $X=j_{n_1}(x_1)j_{n_2}(x_2)....j_{n_k}(x_k)$ with $n \le n_1,.,n_k < \infty$  
by Markov property we have $F_{[n]}XF_{n]}=j_{n}(x)$ for some $x \in \clm^F_0$. We claim further that $\clm^F_{[n}$ is dense in $\clm_{[n}$ in $\sigma-$weak operator topology 
which follows by normal property of complete positive map $\tau$ and $\phi$. For a proof, since $\clm^F_0$ is dense in weak operator topology in $\clm_0$ as $\clm_0$ is a von-Neumann 
algebra of type-I with completely atomic center, we trivially get by normality of the map $j_n$ that $j_n(\clm_0)$ is an element in the weak closer of $\clm^F_{[n}$. Since weak$^*$
closer of $\clm^F_{[n}$ is a von-Neumann algebra, we get $(\clm^{F}_{[n})''=\clm_{[n}$. 

\vsp 
Since $\vee_{n \in \IZ} \clm_{[n} = \clb(\clh)$ by Theorem 2.6, we can use Kaplansky's density theorem to get positive contractive elements $Y_n \in \clm^F_{[-n}$ so that $Y_n \raro |f><f|$ 
as $n \raro \infty$ in strong operator topology. We set $y_n \in \clm_0$ so that $j_{-n}(y_n) = F_{-n]}Y_nF_{-n]}$ and check that $j_{-n}(y_n) \raro F_{-\infty]} |f><f|F_{-\infty]} =|f><f|$ in strong operator 
topology. By our construction and last paragraph, we have now each $y_n \in \clm^F_0$. Let $p'_k$ be the support projection of the positive element $\sum_{1 \le m \le k} y_m$. So in particular 
each $p'_k$ is a finite rank projection and $y_k=y_kp'_k$. $p'_k$ being an increasing projections, it has a limit in strong operator topology say $I-q \in \clm_0$ as $k \uparrow \infty$. If $q \neq 0$, we 
can take a sequence of finite projection $q_k \in \clm_0$ such that $q_k \uparrow q$ and set $p_k=p'_k+q_k$. Thus we have $y_kp_k=y_kp'_k=y_k$ for each $k \ge 1$ and $p_k \uparrow I-q+q=I$ as $k \uparrow \infty$.  

\vsp 
But $j_{-n_k}(y_k)=j_{-n_k}(y_kp_{n_k})=j_{-n_k}(y_k)j_{-n_k}(p_{n_k})$ where we have used $n_k \ge k$ by our construction for $p_{n_k} \ge p_k$ and 
by taking strong operator limit on both side we get $|f><f|=|f><f| |\Omega>< \Omega|=0$. This brings a contradiction. This completes the proof when $\IT$ is 
discrete i.e. $\IZ$. 

\vsp 
For continuous case we note by Kadison-Schwarz inequality $\tau_t(x^*)\tau_t(x) \le \tau_t(x^*x)$ for all $x \in \clm_0$ and so 
$\phi_0(\tau_t(x^*)\tau_t(x)) \le \phi_0(x^*x)$ by invariance of the state. In particular we have $\phi_0(\tau_t(x^*)\tau_t(x)) \le \phi_0(\tau_s(x^*)\tau_s(x))$
for all $t \ge s  \ge 0$. So in order to prove $\phi_0(\tau_t(x^*)\tau_t(x)) \raro |\phi_0(x)|^2$ as $t \raro \infty$, it is enough if we show 
$\phi_0\tau_{tn}(x^*)\tau_{tn}(x)) \raro |\phi(x)|^2$ as $n \raro \infty$ for some fixed $t > 0$. Thus continuous case follows from the discrete case. 
This completes the proof also for $\IT=\IR$.  

\vsp 
That $\clf$ is trivial von-Neumann algebra follows as Kolmogorov's property says now $\phi_0(x^*x)=|\phi_0(x)|^2$ 
for all $x \in \clf$ and $\phi_0$ being faithful $x=\phi_0(x)I$. 

\end{proof}  
\vsp 
\begin{rem}
It is known that unique ground state of $H_{XY}$ model gives an example of a translation invariant pure state $\omega$ on $\clb=\otimes_{\IZ}\!M_d(\IC)$ with a 
type-III factor once restricted to $\clb_0=\otimes_{\IZ_+}\!M_d(\IC)$ and $\clm_0=\pi_{\omega}(\clb_0)''$ i.e. the state $\omega$ is faithful on $\clb_0$. In such 
a case $\tau_t:\clm_0 \raro \clm_0$ is a semi-group of endomorphisms. This in particular shows that Kolmogorov property is not guaranteed by purity property in general. 
Further $\clf=^{\mbox{def}}\{x: \tau_t(x^*)\tau_t(x)=\tau_t(x^*x),\tau_t(x)\tau(x^*)=\tau_t(xx^*):\;t \ge 0 \}$ equal to $\clm_0$. 
Theorem 2.5 says however $\clg$ is trivial and thus modular group plays a non-trivial role in such a situation. It seems no example 
of a translation invariant state $\omega$ on $\clb$ with Kolmogorov property is known  which gives a type-III factor state on $\clb_0$. We avoid giving details here.
\end{rem}

\begin{rem}
In case $\clm_0$ is a commutative von-Neumann algebra both weak Kolmogorov and Kolmogorov's property coincides as modular 
operator is trivial and the argument used in the faulty proof given for Theorem 3.4 in [Mo2] goes through due to triviality of modular operator. On the other hand weak Kolmogorov property of $(\tau_t)$ 
implies in particular $\tilde{\tau}_t(x) \raro \phi_0(x)$ in weak operator topology and so by duality $\tau_t(x) \raro \phi_0(x)$ in weak operator topology as $t \raro \infty$. This shows strong mixing 
is a weaker notion then weak Kolmogorov property and we have following order: 
$$\mbox{Ergodic} << \mbox{Weakly mixing} << \mbox{Strongly mixing} << $$ 
$$\mbox{Weakly Kolmogorov=Purity} << \mbox{Kolmogorov} $$  
In case $\clm_0$ is a type-I factor, strong mixing implies Kolmogorov's property by Theorem 2.8. Further we recall Theorem 2.5 in [Mo2] that when $\IT=\IR$, i.e. continuous time, all above notions 
coincides once $\clm_0$ is a type-I factor.   
\end{rem}

\section{ Pure inductive limit state: }

\vsp
Let $\clb_0$ be a $C^*$ algebra, $(\lambda_t:\;t \ge 0)$ be a semi-group 
of injective endomorphisms and $\psi$ be an invariant state for $(\lambda_t:t \ge 0)$. 
We extend $(\lambda_t)$ to an automorphism on the $C^*$ algebra $\clb_{[-\infty}$ of the 
inductive limit $$ \clb_0 \raro^{\lambda_t} \clb_0 \raro^{\lambda_t} \clb_0 $$
and extend also the state $\psi$ to $\clb_{[-\infty}$ by requiring $(\lambda_t)$ invariance.
Thus there exists a directed set ( i.e. indexed by $\IT$ , by inclusion $\clb_{[-s} \subseteq \clb_{[-t}$ 
if and only if $t \ge s$ ) of C$^*$-subalgebras $\clb_{[t}$ of $\clb_{[-\infty}$ so that the uniform closure of 
$\bigcup_{s \in \IT} \clb_{[s}$ is $\clb_{[-\infty}$. Moreover there exists an isomorphism
$$i_0: \clb_0 \raro \clb_{[0}$$
( we refer [Sa] for general facts on inductive limit of C$^*$-algebras ). 
It is simple to note that $i_t=\lambda_t \circ i_0$ is an isomorphism of $\clb_0$ onto $\clb_{[t}$ and 
$$\psi_{-\infty} i_t = \psi$$ 
on $\clb_0$. Let $(\clh_{\pi},\pi,\Omega)$ be the GNS space associated with $(\clb_{[-\infty},\psi_{[-\infty})$ and $(\lambda_t)$ be 
the unique normal extension to $\pi(\clb_{[-\infty})''$. Thus the vector state $\psi_{\Omega}(X)=<\Omega, X \Omega>$ is 
an invariant state for automorphism $(\lambda_t)$. As $\lambda_t(\clb_{[0}) \subseteq \clb_{[0}$ for all $t \ge 0$, 
$(\pi(\clb_{[0})'',\lambda_t,\;t \ge 0,\psi_{\Omega})$ is a quantum dynamics of endomorphisms. Let $F_{t]}$ be the support 
projection of the normal vector state $\Omega$ in the von-Neumann sub-algebra $\pi(\clb_{[t})''$. $F_{t]} \in \pi(\clb_{[t})'' 
\subseteq \pi(\clb_{[-\infty})''$ is a monotonically decreasing sequence of projections as $t \raro -\infty$. Let the projection 
$F_{-\infty]}$ be the limit. Thus $F_{-\infty]} \ge [\pi(\clb_{[-\infty})'\Omega] \ge |\Omega><\Omega|$. So $F_{-\infty]} = |\Omega><\Omega|$ ensures
that $\psi$ on $\clb_{[-\infty}$ is in particular pure. We aim to investigate when $F_{-\infty]}=|\Omega><\Omega|$.

\vsp
To that end we set von-Neumann algebra $\cln_0=F_{0]}\pi(\clb_{[0})''F_{0]}$ and define family 
$\{ k_t: \cln_0 \raro \pi(\clb_{[-\infty})'',\; t \in \IT \}$ of $*-$homomorphisms by
$$k_t(x) = \lambda_t(F_{0]}xF_{0]}),\;\; x \in \cln_0$$  
It is a routine work to check that $(k_t:t \in \IT)$ is the unique up to isomorphism ( in the
cyclic space of the vector $\Omega$ generated by the von-Neumann algebra $\{k_t(x): t \in \IT, x \in \cln_0 
\}$ ) forward weak Markov process associated with $(\cln_0,\eta_t,\psi_0)$ where 
$\eta_t(x)=F_{0]}\lambda_t(F_{0]}xF_{0]})F_{0]}$ for all $t \ge 0$. It is minimal once restricted to the cyclic space generated 
by the process. Thus $F_{-\infty]} =|\Omega><\Omega|$ when restricted to the cyclic subspace of the process if and only if 
$\psi_0(\eta_t(x)\eta_t(y)) \raro \psi_0(x)\psi_0(y) $ as $t \raro \infty$ for all $x,y \in \cln_0$. 
 
\vsp
\begin{pro} 
Let $G_{0]}$ be the cyclic subspace of the vector $\Omega$ generated by $\pi(\clb_{[0})''$. 

\NI (a) $G_{0]} \in \pi(\clb_{[0})'$ and the map $h: \pi(\clb_{[0})'' \raro G_{0]}\pi(\clb_{[0})''G_{0]}$ defined by
$X \raro G_{0]}XG_{0]}$ is an homomorphism and the range is isomorphic to $\pi_0(\clb_0)''$, where 
$(\clh_{\pi_0},\pi_0)$ is the GNS space associated with $(\clb_0,\psi)$.  

\NI (b) Identifying the range of $h$ with $\pi_0(\clb_0)''$ we have
$$ h \circ \lambda_t(X)=\lambda_t(h(X)) $$
for all $X \in \pi(\clb_{[0})''$ and $t \ge 0$. 

\NI (c) Let $P$ be the support projection of the state $\psi$ in von-Neumann algebra $\pi_0(\clb_0)''$ and 
$\clm_0=P\pi_0(\clb_0)''P$. We set $\tau_t(x)=P\lambda_t(PxP)P$ for all $t \ge 0,\; x \in \clm_0$ and 
$\psi_0(x)=\psi(PxP)$ for $x \in \clm_0$. Then 

\NI (i) $h(F_{0]})=P$ and $h(\cln_0)=\clm_0$;

\NI (ii) $h(\eta_t(x)) = \tau_t(h(x))$ for all $t \ge 0$.   
\end{pro} 

\begin{proof} 
For a proof we refer [Mo2]   
\end{proof} 

\vsp
\begin{thm} 
$\psi_{[-\infty}$ is a pure state if and only if for each $t \ge 0$, there exists contractive elements 
$x_t \in \clm_0$ such that for each $s \ge 0$, $\phi_0(\tau_t(x)x_{s+t}^*x_{s+t}\tau_t(y)) \raro \phi_0(x)\phi_0(y)$ as $t \raro \infty$ for 
$x,y \in \clm_0$. 
\end{thm} 

\begin{proof} 
Let $Q_0$ be the support projection of $\vee_{t \in \IR } \clm_{[t}$. Thus $Q_0 \le F_{-\infty]}  \le F_{t]}$ for all $t \in \IR$. 
For any fix $t \in \IT$ since $k_t(\clm_0)= F_{t]}\pi(\clb_{[t})''F_{t]}$, for any $X \in \pi(\clb_{[t})''$ we have 
$Q_0X\Omega=Q_0F_{t]}XF_{t]}\Omega=Q_0k_t(x)\Omega$ for some $x \in \cln_0$. Hence $Q_0=|\Omega><\Omega|$ if and only if $Q_0=|\Omega><\Omega|$ on the 
cyclic subspace generated by $ \{ k_t(x),\;t \in \IT, x \in \clm_0 \}$. Theorem 2.4 says now that $Q_0=|\Omega><\Omega|$ if for each $s \ge 0$ 
we have $\psi_0(\eta_t(x)x_{s+t}^*x_{s+t}\eta_t(y)) \raro \psi_0(x)\psi_0(y)$ as $t \raro \infty$ for all $x,y \in \cln_0$. 
Since $h$ is a homomorphism and $h \eta_t(x)= \tau_t(h(x))$, we also have $h(\eta_t(x))x_{s+t}^*x_{s+t}\eta_t(y))=\tau_t(h(x))h(x_{s+t}^*)h(x_{s+t})\tau_t(h(x))$. 
Since $\phi_0 \circ h = \psi_0$ we conclude the ``if part'' of the statement identifying $\clm_0$ with $\cln_0$. For the converse we use Kaplansky's density theorem 
to ensure contractive elements $X_t \in i_{-t}(\clb_0)''$ so that $X_t \raro |\Omega><\Omega|$ as $t \raro \infty$ in strong operator 
topology. Now we set $k_{-t}(x_t) = F_{-t]}X_tF_{-t]} \in \clm_{[-t}$ and since $k_{-t}(x_t) \raro |\Omega><\Omega|$ as 
$t \raro \infty$ in strong operator topology we have $k_{-t}(x^*_tx_t) \raro |\Omega><\Omega|$ in weak operator topology. 
This completes the proof once we compute $<\Omega,k_{-s}(x)k_{-t}(x_t^*x_t)k_{-s}(y)\Omega>=\phi_0(\tau_{t-s}(x^*)x^*_tx_t\tau_{t-s}(y))$ for any $t \ge s$ and $x,y \in \clm_0$. 
\end{proof} 

\vsp 
The following theorem includes a proof for Theorem 1.1. 
 
\begin{thm}  
Let $(\clb_0,\lambda_t:t \ge 0,\psi_0)$ be a unital semi-group of injective endomorphisms and $\psi_0$ be a factor state. 
If $\psi_{[-\infty}$ is pure then $\psi_0$ is either type-I or type-III factor state.   
\end{thm}  

\begin{proof} 
We set $P\pi_0(\clb_0)''P=\clm_0$ where $P$ is the support projection of the state $\psi_0$ in the GNS space i.e. $P=[\pi_0(\clb_0)'\Omega]$. 
Thus $\clm_0$ is a factor of same type. $\psi_{[-\infty}$ being a pure state on $\clb_{[-\infty}$, Theorem 2.1 and Theorem 3.2 says that $\psi_{[-\infty}$ is also pure on $\clm_{[-\infty}$. 
Thus the statement is a simple consequence of Theorem 2.5.  
\end{proof}

\bigskip
{\centerline {\bf REFERENCES}}

\begin{itemize}
\bigskip 

\item{[Ac1]} Accardi, L.: The non-commutative Markov property. (Russian) Funkcional. Anal. i Priložen. 9 (1975), no. 1, 1-8.

\item{[Ac2]} Accardi, L.: Non-commutative Markov chains associated to a preassigned evolution: an application to the quantum theory of measurement. 
Adv. in Math. 29 (1978), no. 2, 226-243. 

\item{[AcC]} Accardi, Luigi; Cecchini, Carlo: Conditional expectations in von Neumann algebras and a theorem of Takesaki.
J. Funct. Anal. 45 (1982), no. 2, 245-273. 

\item {[AFL]} Accardi, L, Frigerio, A, Lewis, John T. Quantum stochastic processes. Publ. Res. Inst. Math. Sci. 18 (1982), 
no. 1, 97-133. 

\item {[AM]} Accardi, L., Mohari, A.: Time reflected Markov processes. Infin. Dimens. Anal. Quantum Probab. Relat. Top., vol-2, no-3, 
397-425 (1999).

\item{[AraMa]} Araki, H., Matsui, T.: Ground states of the XY model, Commun. Math. Phys. 101, 213-245 (1985).

\item{[Ar1]} Arveson, W.: Pure $E_0$-semi-groups and absorbing states. Comm. Math. Phys. 187  (1997),  no. 1, 19--43.

\item{[Ar2]} Arveson, W.: Noncommutative Dynamics and $E_0$-semi-groups, Springer Monogr. Math., Springer, New York, 2003.  

\item{[BhP]} Bhat, B.V.R., Parthasarathy, K.R.: Kolmogorov's existence theorem for Markov processes on $C^*$-algebras, Proc.
Indian Acad. Sci. 104,1994, p-253-262.

\item{[BR]} Bratteli, Ola., Robinson, D.W. : Operator algebras and quantum statistical mechanics, I,II, Springer 1981.

\item{[BJP]} Bratteli, Ola; Jorgensen, Palle E. T.; Price, Geoffrey L. Endomorphisms of B(H). Quantization, nonlinear partial differential equations, 
and operator algebra (Cambridge, MA, 1994), 93–138, Proc. Sympos. Pure Math., 59, Amer. Math. Soc., Providence, RI, 1996. 

\item {[De]} Dell'Antonio, G.F.: On the limit of sequences of normal states, Comm. Pure Appl. Math. 20 (1967) 413-429. 

\item {[Di]} Dixmier, J.: Von Neumann Algebras, North-Holland, 1981. 

\item {[EvKa]} Evans, D. E.; Kawahigashi,Y.: Quantum symmetries on operator algebras, 
Oxford University Press. 

\item{[EvL]} Evans, D. E.; Lewis, J. T.: Dilations of irreversible evolutions in algebraic quantum theory.  
Comm. Dublin Inst. Adv. Studies Ser. A No. 24 (1977), v+104 pp.   

\item{[Ha]} Haagerup, U.: The standard form of von-Neumann algebras, Pre-print Copenhagen, also Math. Scand. 37 (1975), no. 2, 271-283. 

\item{[Ma]} Matsui, T.: Private communication (2011).

\item{[Mac]} Mackey, G.W.: Imprimitivity for representations of locally compact gropups I, Proc. Nat. Acad. Sci. U.S.A. 35 (1949), 537-545. 

\item{[Mo1]} Mohari, A.: Markov shift in non-commutative probability, Jour. Func. Anal. 199 (2003) 189-209.  

\item{[Mo2]} Mohari, A.: Pure inductive limit state and Kolmogorov's property, J. Func. Anal. vol 253, no-2, 584-604 (2007)
Elsevier Sciences.

\item{[Mo3]} Mohari, A: Jones index of a completely positive map, Acta Applicandae Mathematicae. Vol 108, Number 3, 665-677 
(2009)

\item{[Mo4]} Mohari, A:  Translation invariant pure state on $\clb=\otimes_{\IZ}M_d(\IC)$ and its split property, arXiv:1310.1886.

\item{[Mo5]} Mohari, A.: A complete weak invariance for Kolmogorov states on $\clb = \dsp{\otimes_{k \in \IZ}}\!M^{(k)}_d(\IC)$, arXiv:1309.7606.

\item {[NeS]} Neshveyev, S.: St\o rmer, E. : Dynamical entropy in operator algebras. Springer-Verlag, Berlin, 2006. 

\item {[Or]} Ornstein, D..: Two Bernoulli shifts with infinite entropy are isomorphic. Advances in Math. 5 1970 339-348 (1970).  

\item {[OP]} Ohya, M., Petz, D.: Quantum entropy and its use, Text and monograph in physics, Springer-Verlag 1995. 

\item{[Pa]} Parry, W.: Topics in Ergodic Theory, Cambridge University Press, 1981. 

\item{[Po1]} Powers, R. T.: Representation of uniformly hyper-finite algebras and their associated von-Neumann rings, Ann. Math. 86 (1967), 138-171. 

\item{[Po2]} Powers, R. T.: An index theory for semi-groups of $*$-endomorphisms of
$\clb(\clh)$ and type II$_1$ factors.  Canad. J. Math. 40 (1988), no. 1, 86-114.

\item {[Sak]} Sakai, S.: C$^*$-algebras and W$^*$-algebras, Springer 1971.  

\item {[Sa]} Sauvageot, Jean-Luc: Markov quantum semi-groups admit co-variant Markov $C^*$-dilations. Comm. Math. Phys. 
106 (1986), no. 1, 91-103.

\item{[Ta]} Takesaki, M. : Theory of Operator algebras II, Springer, 2001.

\item {[Vi]} Vincent-Smith, G. F.: Dilation of a dissipative quantum dynamical system to a quantum Markov process. Proc. 
London Math. Soc. (3) 49 (1984), no. 1, 58-72.

\end{itemize}

\end{document}